\documentclass{amsart}\usepackage{amsmath,amssymb,amsthm}

\newtheorem{theorem}{Theorem}
\newcommand{\bt}{\begin{theorem}}
\newcommand{\et}{\end{theorem}}
\newtheorem{lemma}{Lemma}
\newcommand{\bl}{\begin{lemma}}
\newcommand{\el}{\end{lemma}}
\newtheorem{corollary}{Corollary}
\newcommand{\bc}{\begin{corollary}}
\newcommand{\ec}{\end{corollary}}
\newcommand{\bconj}{\begin{conjecture}}
\newcommand{\econj}{\end{conjecture}}
\newtheorem{problem}{Problem}
\newcommand{\bprob}{\begin{problem}}
\newcommand{\eprob}{\end{problem}}
\newcommand{\beq}{\begin{equation}}
\newcommand{\eeq}{\end{equation}}
\newcommand{\benum}{\begin{enumerate}}
\newcommand{\eenum}{\end{enumerate}}
\newcommand{\N}{\ensuremath{ \mathbf N }}
\newcommand{\Z}{\ensuremath{\mathbf Z}}

\newcommand{\R}{\ensuremath{\mathbf R}}

\newcommand{\mch}{\ensuremath{ \mathcal H}}

\newcommand{\mcr}{\ensuremath{ \mathcal R}}

\newcommand{\mct}{\ensuremath{ \mathcal T}}

\newcommand{\mba}{\ensuremath{ \mathbf a}}

\newcommand{\mbu}{\ensuremath{ \mathbf u}}

\newcommand{\mbx}{\ensuremath{ \mathbf x}}
\newcommand{\mby}{\ensuremath{ \mathbf y}}

\DeclareMathOperator{\card}{\text{card}}

\DeclareMathOperator{\diam}{\text{diam}}

\newcommand{\bmat}{\left(\begin{matrix}}
\newcommand{\emat}{\end{matrix}\right)}
\newcommand{\bsmallmat}{\left(\begin{smallmatrix}}
\newcommand{\esmallmat}{\end{smallmatrix}\right)}
\DeclareMathOperator{\qand}{\quad\text{and}\quad}
\DeclareMathOperator{\qqand}{\qquad\text{and}\qquad}

\title{Diversity, equity, and inclusion for problems  in additive number theory}

\author{Melvyn B.  Nathanson}
\address{ Lehman College (CUNY), Bronx, NY 10468} 
\email{melvyn.nathanson@lehman.cuny.edu}

\subjclass[2000]{11B13, 11B05, 11B75,  11P70, 22D99} 
\keywords{ Sumsets in semigroups,  restricted sumsets, intersections of sumsets, 
distribution of of sumset sizes,  computational complexity, 
additive number theory, combinatorial number theory}

\thanks{This is the text of a talk delivered at the New York Number Theory Seminar 
at the CUNY Graduate Center on February 19, 2026, and at the International Number Theory 
Conference (Alladi70) at the University of Florida on March 18, 2026.}

\begin{document}

\maketitle

\begin{abstract}
 There is a diversity of  problems in additive number theory. 
 Equity requires the consideration of less well known problems, 
 and suggests their inclusion in the additive canon.  
 Of particular interest are problems about the sizes of sumsets  of finite sets of integers 
 and problems about the arithmetical structure of  intersections of sumsets. 
\end{abstract}

\section{Canonical additive problems} 
Here are three of the most important and best known  research topics in additive number theory:
\benum
\item
Partitions and the associated algebraic, analytic, and combinatorial methods,

\item
Waring's problem on sums of $k$th powers and the circle method,

\item
Additive prime number theory: Sums and differences of primes (e.g. the Goldbach and twin prime conjectures) 
and sieve methods.    
\eenum

These are the canonical additive problems, 
which I shall ignore them and focus, instead,  on other  less popular and more ``diverse'' 
questions that have not always been or  have not yet been included in the additive canon.  
Fairness, or equity, suggests that less well-known and not obviously ``important'' 
topics should always be considered candidates for inclusion, because 
only time and not today's ``experts''  decide what is really important.  
 (For surveys  of the traditional additive core, see Andrews, \emph{The Theory of Partitions}~\cite{andr84} 
and Nathanson, \emph{Additive Number Theory: The Classical Bases}~\cite{nath96aa}.)

\section{Sum-product: A problem no one cared about} 
Paul Erd\H os ~\cite{erdo77}, in the book \emph{Number Theory Day, Proceedings of the Conference Held at Rockefeller University, New York 1976}, wrote the following:

\begin{quotation}
 
Section 7. A new extremal problem 

Let $1 \leq a_1 < \cdots < a_n$ be a sequence of integers.  
Denote by $f_r(n)$ the minimal number of distinct integers which are the sum or 
product of exactly $r$ of the $a$'s. 

I conjectured that for every $r$ and $\varepsilon > 0$ if $n > n(\varepsilon,r)$, 
$f_r(n) > n^{r-\varepsilon}$.

If true this seems extremely difficult and seems to require new ideas (unless of course an obvious point is being overlooked).

Szemer\' edi and I observed that it follows from deep results of Freiman that 
\[
\lim_{n\rightarrow \infty} f_2(n)/n = \infty
\]
but even the proof of $f_2(n) > n^{1+\varepsilon}$ seems to present great difficulties.  
\end{quotation}

Several years later, in 1983,  Erd\H os and  Szemer\' edi~\cite{erdo-szem83} 
proved that there exists an $\varepsilon > 0$ 
such that 
 \[
 f_2(n) > n^{1+\varepsilon}.
 \] 
An important result, but did anyone pay attention?  
In 1991, Erd\" os~\cite{erdo91} wrote that  “our paper with Szemer\' edi has nearly been forgotten.''

In 1997, Nathanson~\cite{nath97} computed  the first explicit value for $\varepsilon$ : 
\[
\varepsilon = \frac{1}{31}  \qqand 
 f_2(n) > c n^{32/31} 
 \]
where $c = 0.00028\ldots$.
The explicit exponent $\varepsilon$ was quickly improved. \\
Ford~\cite{ford98} proved 
\[
\varepsilon = \frac{1}{15} - o(1). 
\]
Elekes~\cite{elek97} proved 
\[
\varepsilon = \frac{1}{4} - o(1)
\]
Solymosi~\cite{soly09} proved 
\[
\varepsilon = \frac{1}{3} - o(1)
\]
Bloom~\cite{bloo25} has the current best result: 
\[
\varepsilon = \frac{1}{3}  + \frac{5}{9813} - o(1)
\]
The ``sum-product'' problem has become part of the canon 
and big business in what is now called ``additive combinatorics.''

\section{Freiman's theorem} 

In the 1976 Number Theory Day conference at Rockefeller University, Erd\H os wrote, 
``Szemer\' edi and I observed that it follows from deep results of Freiman that \ldots.'' 

Almost no one had paid any attention to Grigorii Freiman's remarkable work on inverse problems in 
additive number theory, even though he had  published  a book in Russian in 1966 
and the American Mathematical Society (Freiman~\cite{frei73}) 
had published an English translation in 1973. 
(Of course, in Freiman's book and papers, the proof of his big theorem was hard to understand.) 
Much later, various Hungarians worked on the theorem.  Ruzsa~\cite{ruzs94}   
found a beautiful proof in 1994, which was exposed in the  author's text 
\emph{Additive Number Theory: Inverse Theorems and the 
Geometry of Sumsets}~\cite{nath96bb}.  
Tim Gowers~\cite{gowe01} used this to obtain explicit estimates for Szemer\' edi's theorem on 
$k$-term arithmetic progressions 
in dense sets of integers, and Freiman's theorem and its consequences became another 
big business in additive number theory.

\section{A noncanonical problem: Sizes of sumsets} 
For every set $A$ of integers, consider the 2-fold sumset 
\[
2A = \{ a_i + a_j : a_i, a_j \in A \}.
\]
In their 1983 paper on the sum-product conjecture,  
Erd\H os and  Szemer\' edi~\cite{erdo-szem83} wrote: 
\begin{quotation} 
Let $2k-1 \leq t \leq \frac{k^2+k}{2}$.  It is easy to see that one can find a sequence 
of integers $a_1 < \ldots < a_k$ so that there should be exactly $t$ distinct integers 
in the sequence $a_i+a_j, 1 \leq i \leq j \leq k$. 
\end{quotation} 

The lower bound comes from $k$-term arithmetic progressions.  
For example, every integer in the interval of integers  $\{0,1,2,\ldots, 2k-2\}$ is a sum of two integers in the 
arithmetic progression $A = \{0,1,2,\ldots, k-1\}$.  
Moreover, the set of 2-fold sums of a set of $k$ integers 
has size $2k-1$ if and only if the set is an arithmetic progression.  

The upper bound is combinatorial:  It comes from sets $A = \{a_1,\ldots, a_k\}$  of size $k$ 
whose sums $a_i+a_j$ with $a_i \leq a_j$ are distinct, and so 
\[
|2A| =  \binom{k+1}{2} = \frac{k^2+k}{2}.
\]   
Such sets are called \emph{Sidon sets}.  
For example, the set $\{1,3,3^2,\ldots, 3^{k-1}\}$ is a Sidon set.

For all positive integers  $k$, define the \emph{sumset size set} 
\[
\mcr (2,k)  = \left\{ \left| 2A \right|: A \subseteq \Z \text{ and }  |A| = k \right\} . 
 \]
 One can also consider sums of distinct integers. 
Let 
\[
\widehat{2A} = \{ a_i + a_j : a_i, a_j \in A \text{ and } a_i \neq a_j \} 
\]
and  
\[
\widehat{\mcr} (2,k) = \left\{ \left|\widehat{2A} \right|: A \subseteq \Z \text{ and }  |A| = k \right\}.
 \] 

For $x<y$,  define the \emph{interval of integers}  $[x,y] = \{ t \in \Z: x \leq t \leq y\}$. 

\bt [Nathanson~\cite{nath26b}]        
For all integers $k \geq 2$, 
\beq        \label{DEI:elementaryEstimate-2}
\mcr (2,k) =   \left[ 2k-1, \frac{k^2 +k}{2} \right].
\eeq
Moreover, for all $t \in \mcr(2,k)$, there exists a set $A \subseteq \left[0,2^k -1 \right]$ 
such that $|A| = k$ and $\left| 2A\right| = t$.  

For  all   integers $k \geq 2$, 
\[
\widehat{\mcr} (2,k) =  \left[ 2k-3, \frac{k^2 - k}{2} \right].
\]
Moreover, for all $t \in\widehat{\mcr}(2,k)$, there exists a set 
$A \subseteq \left[0,2^{k-2} \right]$
such that $|A| = k$ and $\left|\widehat{2A} \right| = t$.
\et

\bprob
Are there   subexponential or  polynomial bounds on the set $A$?  
Do  there exist  integers $n$ and $c$ 
such that, for all $k \geq 2$ and 
 for all $t \in \mcr (2,k)$, there is a set $A \subseteq \left[0, ck^n\right]$ 
with $|A| = k$ and $\left|2A \right| = t$?
Do  there exist  integers $n$ and $c$ such that, for all $k \geq 2$ and 
 for all $t \in \widehat{\mcr} (2,k) $, there is a set $A \subseteq \left[0, ck^n\right]$ 
with $|A| = k$ and $\left|\widehat{2A} \right| = t$?
\eprob

\bprob
``Most'' sets $A$ of size $k$  are Sidon sets, that is, $|2A| = \binom{k+1}{2}$.
What is the distribution of sumset sizes $|2A|$?
What are the ``popular'' sizes and why are  they popular?

What is the distribution of sumset sizes $\left|\widehat{2A} \right| $? 
What are the ``popular'' sizes and why are  they popular?
\eprob

\section{Range of sumset sizes}
%Let $A$ be a nonempty finite set of integers with  $\min(A) = 0$ and $|A| = k$.
  
The $h$-fold \emph{sumset}  of a set $A$ of integers is the set of all sums of $h$ 
not necessarily distinct elements of $A$:
\begin{align*} 
hA & = \underbrace{A+\cdots + A}_{\text{$h$ summands}} \\
& = \left\{ a_1+\cdots + a_h: a_i \in A \text{ for all } i \in [1,h] \right\}. 
\end{align*} 
The $h$-fold \emph{restricted sumset}  of $A$ 
is the set of all sums of $h$  distinct elements of $A$:
\begin{align*}
\widehat{hA} 
& = \left\{ a_1+\cdots + a_h: a_i \in A \text{ for all } i \in [1,h] \right. \\
& \left .  \qquad  \text{ and } a_i \neq a_j \text{ for all } i \neq j \right\} \\ 
& = \left\{ \sum_{s\in S}s :   S \subseteq A \text{ and } |S| = h \right\}. 
\end{align*}
%where $|S|$ denotes the size of the set $S$. 
The sumset $hA$ and the restricted sumset $\widehat{hA}$ are finite if $A$ is finite.

For all positive integers $h$ and $k$, we define the \emph{sumset size sets} 
\[
\mcr (h,k) = \left\{ \left| hA \right|: A \subseteq \Z \text{ and }  |A| = k \right\}   
 \]
and 
\[
\widehat{\mcr} (h,k) = \left\{ \left|\widehat{hA} \right|: A \subseteq \Z \text{ and }  |A| = k \right\}.
 \] 
In general, one considers the range of sumset sizes of finite subsets 
of an arbitrary additive abelian group or semigroup $X$:
\[
\mcr_{X} (h,k) = \left\{ \left| hA \right|: A \subseteq X \text{ and }  |A| = k \right\} 
 \]
and 
\[
\widehat{\mcr}_X (h,k) = \left\{ \left|\widehat{hA} \right|: A \subseteq X \text{ and }  |A| = k \right\}.
 \] 
 
For simplicity, we consider only $\mcr (h,k)$. 
Analogous to~\eqref{DEI:elementaryEstimate-2}, there is the elementary estimate:
\beq        \label{DEI:elementaryEstimate}
\mcr(h,k) \subseteq  \left[ h(k-1)+1,\binom{h+k-1}{h}\right]. 
\eeq 
The lower bound comes from $k$-term arithmetic progressions and the upper bound from 
sets (called $B_h$-sets)  with unique $h$-fold sums.

\section{Missing sumset sizes for 3-fold sums of sets with 3 elements}
From the elementary estimate~\eqref{DEI:elementaryEstimate}, we obtain 
\[
\mcr(3,3) \subseteq \{7,8,9,10\}.
\]
Extensive computer searches showed 
\[
\{7,9,10\} \subseteq \mcr(3,3) 
\]
but produced no set $A$ such that $|A| = 3$ and $|3A| = 8$.  
Question: Does there exist a set $A$ with $|A| = 3$ and $|3A| = 8$?

Answer: No, but why not?  

\bt[Nathanson~\cite{nath25a}] 
\[
\mcr(3,3) = \{7,9,10\}. 
\]
\et

Isaac Rajagopal obtained the following beautiful result for 3-fold sumsets:

\bt[Rajagopal~\cite{raja25}] 
For all $k \geq 3$, 
\[
\mcr(3,k)= \{3k-2\} \cup \left[ 3k, \binom{k+2}{3}\right].
\]
\et

\bprob
For all positive integers $k$, 
\[
\mcr(4,k) \subseteq \left[4k-3,\binom{k+3}{4} \right].
\]
Compute $\mcr(4,k)$ for all $k$.   
\eprob

The first missing sumset size result was the following: 

\bt[Nathanson~\cite{nath25a}] 
For all $h \geq 3$ and $k \geq 3$, 
\[
 hk -h + 1 = \min\left( \mcr (h,k) \right) 
\]
but
\[
hk -h + 2 \notin \mcr (h,k). 
\]
\et

This was improved by Vincent Schinina. 

\bt[Schinina~\cite{schi25}]
For all $h \geq 3$ and $k \geq 3$, 
\[
 hk -h + 1 \in  \mcr (h,k)  \qand hk  \in  \mcr (h,k)  \qand 
\]
but
\[
[hk -h + 2, hk-1] \cap \mcr (h,k) = \emptyset.  
\]
\et

Thus, the set   $\mcr(h,k)$ 
is \textbf{never} an interval of consecutive integers for $h \geq 3$ and $k \geq 3$. 
This suggests the question:

\bprob 
What are other ``missing'' sumset sizes?
\eprob

The following result is fundamental.

\bt[Rajagopal~\cite{raja25}]           \label{DEI:theorem:Rajagopal-2}
Let 
\[
m = h(k-1)+1
\]
and 
\[
\Delta_{h,k} = \bigcup_{\ell=1}^{\min(h,k)-3} [m + \ell h + 1, m + \ell h + (h-2-\ell)].
\]
For all $h$ and $k$,
\[
\mcr(h,k) \cap \Delta_{h,k} = \emptyset.
\]
For all $h$, there exists $k_h$ such that, for all $k > k_h$, 
\[
\mcr(h,k)= \left[ h(k-1)+1,\binom{h+k-1}{h}\right] \setminus  \Delta_{h,k}. 
\]
\et

\bprob
Compute $k_h$ in Theorem~\ref{DEI:theorem:Rajagopal-2}. 
\eprob

\section{Computational complexity}

Let $h \geq 2$ and $k \geq 2$.  Because the set   $\mcr(h,k)$ is finite, 
there is an integer $N$ such that,  if $t \in \mcr (h,k)$, 
then there exists $A \subseteq [0,N]$ 
with $|A| = k$ and $|hA| = t$.

\bprob
Compute an integer $N$ with this property. 
Compute the smallest integer $N(h,k)$ with this property.
\eprob

Noah Kravitz and Nathanson obtained an upper bound. 

\bt[Nathanson~\cite{nath26a}]
There is the upper bound 
\[
N(h,k) < 4(4h)^{k-1}.
\]
\et

For fixed $k$, this is polynomial in $h$.  For fixed $h$, this is exponential in $k$.  
Thus,  $N(h,4) < 256h^3$ and  $N(4,k) <  4^{2k-1}$.

\bprob
Improve these bounds. 
\eprob

\bprob
If $A$ is a set of size $k$, then the sumset $hA$ and the size 
of the sumset $hA$ can be computed in polynomial time, 
and so the problem of determining the sumset size set $\mcr(h,k)$ 
is in the complexity 
class NP for all $h$ and $k$.  

Prove (or disprove): This problem  is neither  NP-hard nor NP-complete.  
\eprob

\section{Sumset size races}

Here is another class of sumset size problems in additive number theory. 

Let $A_1$, $A_2$, \ldots, $A_n$ be $n$ sets of $k$ integers.  
For every positive integer $h$, consider the $n$-tuple of $h$-fold sumset sizes
\[
(|hA_1|, |hA_2|, \ldots, |hA_n|).
\]
If the sumset sizes are pairwise distinct, we can associate to this $n$-tuple 
a permutation of the integers $1,2,\ldots, n$ which describes the relative order 
of the sequence of sumset sizes.  
If the sumset sizes are not pairwise distinct, we can still associate to this $n$-tuple 
an $n$-tuple of positive integers that  describes the relative order 
of the sequence of sumset sizes.

\bprob
Given $H$ not necessarily distinct $n$-tuples of positive integers 
$\tau_1, \tau_2, \ldots, \tau_H$,  do there exist finite sets 
$A_1, A_2, \ldots,  A_n$ of size $k$ such that, for all $h = 1,2,\ldots, H$, 
the relative   order of the sumset size $n$-tuple 
\[
(|hA_1|, |hA_2|, \ldots, |hA_n|).
\]
is the same as the relative order of the $h$-th integer $n$-tuple $\tau_h$?
\eprob

One can ask the same question about sumsets in arbitrary additive abelian groups.
 
By theorems of Nathanson~\cite{nath72} (for integers) and Khovanskii~\cite{khov95} (for groups), 
for an arbitrary infinite sequence $\tau_1, \tau_2, \ldots, $ of  $n$-tuples of positive integers , 
there do not exist finite sets 
$A_1, A_2, \ldots,  A_n$ of size $k$ such that the relative order of the $h$-fold sumset size $n$-tuple 
\[
(|hA_1|, |hA_2|, \ldots, |hA_n|).
\]
is the same as the relative order of the $h$-th integer $n$-tuple $\tau_h$ for all $h = 1,2,\ldots$.

Following work of Paul P\' eringuey and Anne de Roton~\cite{peri-roto25}, 
Noah Kravitz proved the following:

\bt[Kravitz~\cite{krav25}]          \label{DEI:theorem:Kravitz2}
Let  $\tau_1, \tau_2, \ldots, \tau_H, \tau_{\infty}$ be  $n-$tuples of positive integers.  
There exist  finite sets $A_1, \ldots, A_n$ such that the sumset size $n$-tuple 
\[
(|hA_1|, |hA_2|, \ldots, |hA_n|) 
\]
has the same relative order as the $n$-tuple $\tau_h$ for all $h=1,\ldots, H$, 
and the sumset size $n$-tuple has the same relative order 
as the $n$-tuple $\tau_{\infty}$ for all $h > H$. 
\et

\bc
If $\varepsilon_h = 1, 0,$ or $-1$ for all $h = 1,\ldots, H$, then 
there exist finite sets $A$ and $B$ such that  
\[
|hA| - |hB|  
\]
has the same sign as $\varepsilon_h$ for all $h = 1,\ldots, H$. 
\ec

There is the following refinement by Jacob Fox, Noah Kravitz, and Shengtong Zhang.

\bt[Fox,  Kravitz, Zhang~\cite{fox-krav-zhan25}]
For every sequence
$t_1,\ldots, t_H$ of integers, 
there exist finite sets $A$ and $B$ such that  
\[
|hA| - |hB| = t_h
\]
for all $h=1,\ldots, H$. 
\et

This is true not just for finite sets of integers but for finite subsets of infinite or 
of  sufficiently large finite abelian groups.

\section{Sumset intersection limits} 
Let $\N = \{1,2,3,\ldots\}$ be the set of positive integers 
and $\N_0 = \{0,1,2,3,\ldots\}$ be the set of nonnegative  integers. 
A decreasing sequence of sets $(A_q)_{q=1}^{\infty}$ is 
\emph{strictly decreasing} if $A_q \neq A_{q+1}$ for all $q \in \N$ 
and \emph{asymptotically strictly decreasing} if $A_q \neq A_{q+1}$ 
for infinitely many $q \in \N$. 
A decreasing sequence of sets $(A_q)_{q=1}^{\infty}$ that is not asymptotically 
strictly decreasing is eventually constant. 

Let $(A_q)_{q=1}^{\infty}$ be a decreasing sequence of nonempty subsets 
of an additive abelian semigroup $X$.  
Then  $(hA_q)_{q=1}^{\infty}$ is a decreasing sequence of sumsets.  
We define 
\[
A = \lim_{q\rightarrow \infty} A_q =   \bigcap_{q=1}^{\infty} A_q. 
\] 
For all $q \in \N$, we have $A \subseteq A_q$ and so $hA \subseteq hA_q$.  
Therefore, 
\[
hA \subseteq \bigcap_{q=1}^{\infty} hA_q.
\]
Equivalently, 
\[
 h   \lim_{q\rightarrow \infty} A_q \subseteq \lim_{q\rightarrow \infty} hA_q   
 \]
for every sequence $(A_q)_{q=1}^{\infty}$. 
It is natural to ask:  For what sequences $(A_q)_{q=1}^{\infty}$ do we have 
\[
hA  =  h   \lim_{q\rightarrow \infty} A_q = \lim_{q\rightarrow \infty} hA_q 
 = \bigcap_{q=1}^{\infty} hA_q?  
 \]

For $x \in X$, the \emph{representation function} $r_{A,h}(x)$ counts 
the number of $h$-tuples 
$(a_1,\ldots, a_h) \in A \times \cdots \times A$ such that
$x = a_1 + \cdots + a_h$.  The representation function determines the sumset: 
 $hA = \{x\in X: r_{A,h}(x) > 0\}$. 
 
 In the additive semigroup $\N_0$, we have $r_{\N,h}(x) < \infty$ for all $x \in \N_0$ and $h \geq 2$.
  In the additive semigroup \Z, we have $r_{\Z,h}(x) = \infty$ for all $x \in \Z$ and $h \geq 2$.

\bt[Nathanson~\cite{nath26c,nath26e} ]            \label{intersect:theorem:finiteReps}
Let $h \in \N$.  Let $X $ be an additive abelian semigroup such that $r_{X ,h}(x) < \infty$ 
for all $x \in X $.  Let $A$ be a subset of $X$.  Then 
\[
hA = \bigcap_{q=1}^{\infty} hA_q 
\] 
for every asymptotically strictly decreasing sequence $(A_q)_{q=1}^{\infty}$ 
of subsets of $X$ such that $A = \bigcap_{q=1}^{\infty} A_q$. 

If $r_{X ,h}(x) < \infty$ for all $x \in X$ and all $h \geq 2$, then 
\[
hA = \bigcap_{q=1}^{\infty} hA_q 
\]
for all $h \geq 2$ and for every asymptotically strictly decreasing sequence $(A_q)_{q=1}^{\infty}$ 
of subsets of $X$ such that $A = \bigcap_{q=1}^{\infty} A_q$.  
\et

Let $(A_q)_{q=1}^{\infty}$ be a strictly decreasing sequence of sets 
in an additive abelian semigroup $X$ and let 
$A = \bigcap_{q=1}^{\infty}A_q$.  
Consider the set 
\begin{align*}
\mch(A_q) 
& = \{h \in \N: hA = \bigcap_{q=1}^{\infty} hA_q\} \\ 
& = \left\{ h \in \N:  h   \lim_{q\rightarrow \infty} A_q = \lim_{q\rightarrow \infty} hA_q \right\}. 
\end{align*} 
By Theorem~\ref{intersect:theorem:finiteReps}, in the additive semgroup $\N_0$ 
of nonnegative integers, we have $\mch(A_q) = \N$ 
for every sequence $(A_q)_{q=1}^{\infty}$.  This is not necessarily true for other semigroups.   
It is not  even true in the group \Z\ of integers that $\mch(A_q) = \N$ 
for every asymptotically decreasing sequence  $(A_q)_{q=1}^{\infty}$. 
For example, there is the following result.

\bt[Marques-Nathanson~\cite{marq-nath26}]
For every $h_0 \geq 2$,  
there exist asymptotically strictly decreasing sequences  $(A_q)_{q=1}^{\infty}$ of sets of integers with 
$\{1,\ldots, h_0-1\} \subseteq \mch(A_q)$ but $h_0 \notin \mch(A_q)$,  
and also that there exist sequences  $(A_q)_{q=1}^{\infty}$ with 
$\{1,  h_0 \} \subseteq \mch(A_q)$ but $\{2,3, \ldots, h_0-1\} \cap \mch(A_q) = \emptyset$.
\et

\bprob 
Describe the sets $\mch(A_q)$ for asymptotically strictly decreasing sequences  
$(A_q)_{q=1}^{\infty}$ of sets of integers. 
\eprob

\bprob 
Compute the set 
\[
\mch^*(0) 
 = \left\{ \mch(A_q) : \bigcap_{q=1}^{\infty}A_q = \{ 0\} \right\}.
\]   
\eprob

\section{The moral of the story}

The problems in this paper may interest only a minority of number 
theorists and, perhaps, only a minority of additive number theorists 
and additive combinatorialists, but, as the history of the  sum-product conjecture 
and Freiman's theorem suggests, inclusion of a diversity of problems into the 
additive number theory discussion may be not only equitable but even beneficial  
to the larger mathematical community.

It suggests one of Aesop's fables~\cite{aesop}. 
Aesop lived in ancient Greece between 620 and 564 BCE.  
 Herodotus, in the fifth century BCE, mentioned him 
 in his  \emph{Histories} (Book II, Chapter 134), where he  referred to him 
 as ``Aesop the maker of fables'' and writes that  he was a slave of Iadmon, a Samian.  
Aesop's  fables were transmitted orally and not collected and written down until 
 three centuries after Aesop's death. 
  Stories have been added continuously and are still being added, 
so it is not clear which of Aesop's fables were recited by Aesop.  

One of the oldest and most famous of Aesop's fables is variously known as 
``The Old Man and His Sons'' 
or ``The Bundle of Sticks''.  
It is numbered 53 in the Perry Index~\cite{perr52}.  

\begin{quotation}
There was an old man who had many sons.  The sons were always fighting with each other, 
and the father kept trying unsuccessfully to get them to stop.  
One day he collected a large number of sticks, 
and, with a rope, tied them together into a thick bundle.  
He asked his sons to break the bundle.  One by one they tried and failed.  
Then the father cut the rope and gave his sons the individual sticks, 
which they broke with ease.  ``Working together, you can succeed,'' the father told them.  
\end{quotation}

The moral of the fable is: In unity there is strength.  

This is a good way to think about  the diversity of problems in additive number theory.

\end{document}